\documentclass[11pt]{article}
\usepackage{amsfonts}
\usepackage{latexsym}
\usepackage{amsmath}
\usepackage{amssymb}
\usepackage[authoryear,round]{natbib}
\newcommand{\qed}{\hfill \ensuremath{\Box}} 
\newcommand{\btheta}{\bar{\theta}}
\newcommand{\E}{\mathrm{E}}
\newcommand{\dd}{\mathrm{d}}

\newcommand{\epsN}{{\varepsilon_N}}
\newcommand{\eps}{{\varepsilon}}
\newcommand{\bpi}{\bar{\pi}}
\newcommand{\bpihal}{\bar{\pi}_{\hat{\alpha}}}
\newcommand{\bpial}{\bar{\pi}_{\alpha}}
\newcommand{\hal}{\hat{\alpha}}

\begin{document}
\begin{center}
{\Large Asymptotically minimax Bayesian predictive densities for multinomial models} \\

\vspace{20pt}
{\large Fumiyasu Komaki} \\
{\small Department of Mathematical Informatics} \\
{\small Graduate School of Information Science and Technology, the University of Tokyo} \\
{\small 7-3-1 Hongo, Bunkyo-ku, Tokyo 113-8656, JAPAN} \\
\end{center}

\vspace{10pt}

Dirichlet prior, Jeffreys prior, Kullback-Leibler divergence, latent information prior, reference prior

\vspace{10pt}

\begin{center}
{\bf Summary}
\end{center}

One-step ahead prediction for the multinomial model is considered.
The performance of a predictive density is evaluated by
the average Kullback-Leibler divergence from the true density to the predictive density.
Asymptotic approximations of risk functions of Bayesian predictive densities
based on Dirichlet priors are obtained.
It is shown that a Bayesian predictive density based on a specific Dirichlet prior
is asymptotically minimax.
The asymptotically minimax prior is different from known objective priors such as
the Jeffreys prior or the uniform prior.

\vspace{0.5cm}

\section{Introduction}

We consider one step ahead prediction for the multinomial model.
Suppose that we observe a random variable $x = (x_1, x_2, \ldots, x_{k-1})$ distributed according to the
multinomial distribution
\begin{align*}
p(x|\theta) = &
\binom{N}{x_1, x_2 \dotsb, x_k} \theta_1^{x_1} \theta_2^{x_2}  \dotsb  \theta_k^{x_k}
\end{align*}
where $x_k := N-\sum^{k-1}_{i=1} x_i$, $\theta = (\theta_1,\ldots,\theta_{k-1})$, $\theta_k := 1 - \sum^{k-1}_{i=1} \theta_i$,
and
\[
 \binom{N}{x_1, x_2, \ldots, x_k} := \frac{N!}{x_1! x_2! \cdots x_k!}.
\]
The parameter space is
\[
\Delta
:= \{\theta = (\theta_1,\theta_2,\ldots,\theta_{k-1}) \mid  \theta_i \geq 0 ~ (i = 1,\ldots,k),~ \theta_k := 1-\sum_{i=1}^{k-1} \theta_i \}.
\]

The objective is to predict $y$ distributed according to the the multinomial distribution
\[
p(y|\theta) = \theta^{y_1}_1 \theta^{y_2}_2 \dotsb \theta_k^{y_k}
\]
with index $1$,
where $y = (y_1,\ldots,y_{k-1})$ and $y_k := 1-\sum_{i=1}^{k-1} y_i$,
by using a predictive density $q(y;x)$.

The performance of a predictive density $q(y ; x)$ is evaluated by the risk function
\begin{align}
R(\theta, q(y;x))=& \sum_y \sum_x p (x, y | \theta) \log \frac{p(y|\theta)}{q (y ; x)},
\label{averageKL}
\end{align}
which is the average Kullback-Leibler divergence from the true density $p(y|\theta)$ to the predictive density $q(y;x)$.

When a Dirichlet prior
\begin{align}
\pi_a (\theta) \dd  \theta_1 \dotsb \dd \theta_{k-1}
= \frac{\Gamma (A)}{\Gamma (a_1) \dotsb \Gamma(a_k)}
\theta^{a_1 -1}_1 \dotsb \theta^{a_{k} -1}_{k}
\dd  \theta_1 \dotsb \dd \theta_{k-1},
\label{dirichlet}
\end{align}
where $A := \sum_{i=1}^k a_i$, $a = (a_1,\ldots,a_k)$ and $a_i > 0$ for every $i$,
is adopted, the posterior density
and the Bayesian predictive density are given by
\begin{align*}
p_{\pi_a} & (\theta | x) \dd  \theta_1 \dotsb \dd  \theta_{k-1}
= \frac{\Gamma (N+A)}
{\Gamma (x_1 + a_1) \dotsb \Gamma (x_k + a_k)} \theta_1^{x_1+a_1-1} \dotsb 
\theta_{k}^{x_k+a_k-1}
\dd \theta_1 \dotsb \dd \theta_{k-1},
\end{align*}
and
\begin{align*}
p_{\pi_a} & (y \mid x)
= \int p(y | \theta) p_{\pi_a} (\theta | x)
 \dd \theta_1 \dotsb \dd \theta_{k-1}
= \frac{B(x_1+y_1+a_1,\ldots,x_k+y_k+a_k)}{B(x_1+a_1,\ldots,x_k+a_k)}
= \frac{\sum\limits_{i=1}^k (x_i + a_i)y_i}{N+A},
\end{align*}
respectively, where
\[
B(x_1,\ldots,x_k) := \frac{\Gamma (x_1) \dotsb \Gamma (x_k) }{\Gamma (\sum\limits_{i=1}^k x_i)}.
\]
We define
\begin{align*}
\bpi_\alpha (\theta) \dd  \theta_1 \dotsb \dd \theta_{k-1}
= \frac{\Gamma (k \alpha)}{\{\Gamma (\alpha)\}^k}
\theta^{\alpha -1}_1 \dotsb \theta^{\alpha-1}_{k}
\dd \theta_1 \dotsb \dd \theta_{k-1},
\end{align*}
which is $\pi_a$ with $a_1 = \cdots a_k = \alpha$.

In the present paper, we consider the asymptotics as the sample size $N$ goes to infinity,
and construct a Bayesian predictive density based on a Dirichlet prior that is asymptotically minimax in the sense described below.
It is known that a minimax predictive density for one step ahead prediction for the multinomial model
can be constructed
by using a latent information prior defined as a prior
maximizing the conditional mutual information between $y$ and $\theta$ given $x$; see Komaki (2011).
However, the explicit form of such a prior is difficult to obtain, and we need to develop asymptotic methods.

We consider a sequence of parameter subspaces
\[
\Delta_{\varepsilon_N}
:= \{\theta = (\theta_1,\theta_2,\ldots,\theta_{k-1}) \mid  \theta_i \geq \varepsilon_N ~ (i = 1,\ldots,k),
~\theta_k := 1 - \sum_{i=1}^{k-1} \theta_i \},
\]
where $\{\varepsilon_N\}$ is a decreasing sequence of real numbers
such that
$\lim\limits_{N \rightarrow \infty} \varepsilon_N = 0$ and $0 < \varepsilon_N < 1/k$ for every $N$,
to avoid singularity problems concerning the boundary of the original parameter space $\Delta$.
Then,
$\Delta_{\varepsilon_N} \subset \Delta_{\varepsilon_{N+1}}$,
$\lim\limits_{n \rightarrow \infty} \Delta_{\varepsilon_N} = \Delta$,
and $\theta_i \in [\varepsilon_N, 1- (k-1)\varepsilon_N]$.
Increasing sequences of parameter subspaces converging to the original parameter space
are often considered to construct asymptotic objective priors; see e.\,g.~\citet{BB89}, \citet{CB94}, and \cite{Bernardo05}.

Let $\pi^{(N)}_*$ be a prior on $\Delta_{\varepsilon_N}$ such that the corresponding Bayesian predictive density
$p_{\pi^{(N)}_*}(y \mid x)$ is minimax with respect to the parameter space $\Delta_{\epsN}$.
Thus,
\[
\sup_{\theta \in \Delta_{\varepsilon_N}} R (\theta, p_{\pi^{(N)}_*}(y \mid x))
= \inf\limits_{q} \sup\limits_{\theta \in \Delta_{\epsN}} R(\theta, q(y ; x)).
\]
The existence of such a prior is guaranteed by Theorem 2 in Komaki (2011),
since $p_\pi(x) > 0$ for every $x$ if $\pi \in \mathcal{P}(\Delta_{\epsN})$.
Here, $\mathcal{P}(\Delta_{\epsN})$ is the set of all probability measures on $\Delta_\epsN$.

We show that the Bayesian predictive density based on a Dirichlet prior $\bpihal$
with $\hal := 1+1/\sqrt{6}$ is asymptotically minimax in the sense that
\begin{align}
\biggl| 
\sup_{\theta \in \Delta_{\varepsilon_N}} R (\theta, p_{\pi^{(N)}_*}(y \mid x))
- \sup_{\theta \in \Delta_{\varepsilon_N}} R(\theta, p_{\bpi_{\hat{\alpha}}}(y \mid x)) \biggr|
=& \mathrm{o}(N^{-2})
\label{1-3-1}
\end{align}
if $\{\epsN\}$ satisfies appropriate conditions.

For example, when the model is binomial $(k=2)$, the minimax prior is
$\theta^{1/\sqrt{6}} (1-\theta)^{1/\sqrt{6}} / B(1+1/\sqrt{6},1+1/\sqrt{6})$
and is different from the Jeffreys prior
$\theta^{-1/2} (1-\theta)^{-1/2} / B(1/2,1/2)$ or the uniform prior.

Although the multinomial model is relatively simple,
the results in the present paper could be a prototype for further development of theories on other models.

Closely related but essentially different prediction problems have been extensively studied
in the framework of reference prior and Bayes coding;
see e.\,g.~\citet{IH73}, \citet{Bernardo79}, \citet{CB94}, and \citet{Bernardo05}.
In this setting, the objective is to predict large amount of future observables without using data at hand.
Roughly speaking, the Jeffreys prior is asymptotically minimax under suitable regularity conditions.

In contrast, we consider here one step ahead prediction by using $N$ observed data at hand and consider the asymptotics
as $N$ goes to infinity.
The priors attaining minimax prediction in these two settings are quite different; see
\citet{Komaki04} and \citet{Komaki11}
for discussion on the relation between the two settings,
and see \citet{Clarke07} for various related approaches.

In Section 2, we obtain an asymptotic approximation of risk functions of Bayesian predictive densities
based on Dirichlet priors.
The approximation is uniform on $\Delta_\epsN$.
In Section 3, we prove that
the Bayesian predictive density based on the Dirichlet prior $\bpi_{\hal}$ with $\hal := 1 + 1/\sqrt{6}$
is asymptotically minimax if $\{\epsN\}$ satisfies appropriate conditions.
In Section 4, some discussions are given.

\section{Asymptotic evaluation of the risk function}

In this section, we obtain an asymptotic approximation, which is uniform for $\theta \in \Delta_\epsN$,
of the risk functions of Bayesian predictive densities based on Dirichlet priors.

The risk function \eqref{averageKL} of $p_{\pi_a}(y|x)$ based on $\pi_a$ defined by \eqref{dirichlet} is given by
\begin{align}
R & (\theta, p_{\pi_a}(y \mid x)))
= \sum^k_{i=1} \theta_i \sum^N_{x_i = 0} \binom{N}{x_i} \theta^{x_i}_i (1- \theta_i)^{N-x_i}
\log \left(\frac{\theta_i}{\displaystyle \frac{x_i + a_i}{N+A}} \right) \notag \\
=& \sum_i \theta_i \sum^N_{x_i = 0} \binom{N}{x_i} \theta^{x_i}_i (1-\theta_i)^{N-x_i}
\left\{-\log \left( \frac{N \theta_i+a_i}{N \theta_i + A \theta_i} \right)
- \log \left(\frac{x_i+a_i}{N\theta+a_i}-1+1 \right) \right\} \notag \\
=& - \sum_i \theta_i \log \left(1 + s_i \right)
- \sum_i \theta_i \sum^N_{x_i = 0} \binom{N}{x_i} \theta_i^{x_i} (1-\theta_i)^{N-x_i}
\log \left(w_i + 1 \right),
\label{60-1}
\end{align}
where 
\[
s_i :=\frac{a_i - A \theta_i}{N\theta_i+A\theta_i}
\text{~~~~and~~~~}
w_i := \frac{x_i + a_i}{N\theta_i + a_i} -1 = \frac{x_i - N \theta_i}{N\theta_i + a_i}.
\]
Here, $w_i$ $(i=1,2,\ldots,k)$ are random variables with $\mathrm{E}_{\theta_i} (w_i)  = 0$,
and $\sum_{i=1}^k \theta_i s_i = 0$.

If we fix a true parameter value $\theta$ satisfying $\theta_i \in (0,1)$ for all $i = 1,\ldots,k$,
then it is easy to verify that
\[
 R (\theta, p_{\pi_a}(y \mid x)) = \frac{k-1}{2N} + \mathrm{O}(N^{-2}).
\]
A higher order pointwise approximation of the risk function has been studied; see \citet{Komaki96}.

Here, instead of the pointwise approximation,
we obtain an asymptotic approximation that is uniform for $\theta \in \Delta_\epsN$.

\vspace{0.5cm}

\noindent
Theorem 1. 
Let $p_{\pi_a}(y \mid x)$ be a Bayesian predictive density
based on a Dirichlet prior $\pi_a$ defined by \eqref{dirichlet}.
Suppose that $\{\epsN\}$ be a decreasing sequence of real numbers such that
$\lim\limits_{N \rightarrow \infty} \epsN =0$,
$\lim\limits_{N \rightarrow \infty} N \epsN = \infty$, and $0 < \epsN < 1/k$ for every $N$.
Then, the risk function $R(\theta,p_{\pi_a}(y \mid x))$ satisfies
\begin{align}
\sup_{\theta \in \Delta_\epsN} & \Biggl| R(\theta,p_{\pi_a}(y \mid x))
- \frac{k-1}{2N}
-
\frac{1}{N^2}
\biggl\{\sum_{i=1}^k \frac{1}{12 \theta_i} \left(6a_i^2 - 12a_i + 5 \right)
- \frac{1}{2} A^2 + A - \frac{1}{2} k + \frac{1}{12} \biggr\} \notag \\
&- \frac{1}{N^3} \biggl\{\sum_{i=1}^k \frac{1}{12 \theta_i^2} \left(-4a_i^3 + 18a_i^2 -24a_i + 9 \right)
+ \sum_{i=1}^k \frac{1}{4 \theta_i} \left(-6 a_i^2 + 12 a_i - 5 \right)
+ \frac{1}{3} A^3 - A + \frac{1}{2} k
\biggr\} \notag \\
&- \frac{1}{N^4} \biggl\{\sum_{i=1}^k \frac{1}{120 \theta_i^3} \left(30 a_i^4 -240 a_i^3 + 660 a_i^2 - 720 a_i + 251 \right)
+  \sum_{i=1}^k \frac{1}{2 \theta_i^2} \left(4 a_i^3 - 18 a_i^2 + 24 a_i - 9 \right) \notag \\
&+  \sum_{i=1}^k \frac{1}{12 \theta_i} \left(42 a_i^2 -84 a_i + 35 \right)
- \frac{1}{4} A^4 + A -\frac{1}{2} k - \frac{1}{120}
\biggr\} \Biggr| = \mathrm{O}(N^{-5} \epsN^{-4}).
\label{theorem1}
\end{align}
\qed

\vspace{0.5cm}

The proof is given at the end of this section.

From Theorem 1, we obtain the following corollaries.

\vspace{0.5cm}

\noindent
Corollary 1.
Suppose that $\{\epsN\}$ be a decreasing sequence of real numbers such that
$\lim\limits_{N \rightarrow \infty} \epsN = 0$,
$\lim\limits_{N \rightarrow \infty} N^{\frac{3}{4}} \epsN = \infty$,
and
$0 < \epsN < 1/k$ for every $N$.
Then,
\begin{align}
&\sup_{\theta \in \Delta_{\epsN}} \Biggl|
 R(\theta,p_{\pi_a}) - \frac{k-1}{2N}
- \frac{1}{N^2}
\biggl\{\sum_{i=1}^k \frac{1}{12 \theta_i} \left(6a_i^2 - 12a_i + 5 \right)
- \frac{1}{2} A^2 + A - \frac{1}{2}k + \frac{1}{12} \biggr\} \notag \\
&~~ - \frac{1}{N^3} \sum_{i=1}^k \frac{1}{12 \theta_i^2} \left(- 4a_i^3 + 18a_i^2 -24a_i + 9 \right)
- \frac{1}{N^4} \sum_{i=1}^k \frac{1}{120 \theta_i^3} \left(30 a_i^4 -240 a_i^3 + 660 a_i^2 - 720 a_i + 251 \right) \Biggr| \notag \\
&= \mathrm{o} (N^{-2}).
\label{col1}
\end{align}
\qed

\vspace{0.5cm}

\noindent
Proof. Since $\lim\limits_{N \rightarrow \infty} N^{\frac{3}{4}} \epsN = \infty$, $|N^{-3} \theta_i^{-1}| \leq N^{-3} \eps_N^{-1} = \mathrm{o}(N^{-2})$,
$|N^{-4} \theta_i^{-1}| \leq N^{-4} \eps_N^{-1} = \mathrm{o}(N^{-3})$,
$|N^{-4} \theta_i^{-2}| \leq N^{-4} \eps_N^{-2} = \mathrm{o}(N^{-2})$,
and
$N^{-5} \eps_N^{-4} = N^{-2} (N^{-3/4} \eps_N^{-1})^{4} = \mathrm{o}(N^{-2})$,
we obtain \eqref{col1} from Theorem 1.
\qed

\vspace{0.5cm}

\noindent
Corollary 2.
Suppose that $\{\epsN\}$ be a decreasing sequence of real numbers such that
$\lim\limits_{N \rightarrow \infty} \epsN = 0$,
$\lim\limits_{N \rightarrow \infty} N^{\frac{3}{4}} \epsN = \infty$,
and
$0 < \epsN < 1/k$ for every $N$.
Then, the risk function of the Bayesian predictive density based on a Dirichlet prior $\bpihal(\theta)$,
where $\hal := 1 + 1/\sqrt{6}$, satisfies
\begin{align}
\sup_{\theta \in \Delta_{\epsN}} &
\Biggl| R(\theta,p_{\bpihal}(y \mid x)) - \frac{k-1}{2N}
+ \frac{1}{N^2} \frac{k-1}{12} \biggl\{ 1 + (7+2\sqrt{6}) k \biggr\} \notag \\
& + \frac{1}{N^3} \frac{1}{18\sqrt{6}} \sum_{i=1}^k \frac{1}{\theta_i^2}
+ \frac{1}{N^4} \Bigl( \frac{1}{6 \sqrt{6}} - \frac{11}{720} \Bigr) \sum_{i=1}^k \frac{1}{\theta_i^3} \Biggr|
= \mathrm{o}(N^{-2})
\label{ahat-risk}
\end{align}
and
\begin{align}
\sup_{\theta \in \Delta_{\epsN}} &
\Biggl( R(\theta,p_{\bpihal}(y \mid x)) - \frac{k-1}{2N}
+ \frac{1}{N^2} \frac{k-1}{12} \biggl\{ 1 + (7+2\sqrt{6}) k \biggr\} \Biggr)
= \mathrm{o}(N^{-2}).
\label{ahat-suprisk}
\end{align}
\qed

\vspace{0.5cm}

\noindent
Proof.
We have \eqref{ahat-risk} from Corollary 1
because
$\displaystyle 6 \hal^2 - 12 \hal + 5 = 0$, $\displaystyle - 4 \hal^3 + 18 \hal^2 -24 \hal + 9 = - \sqrt{6}/9$,
$\displaystyle 30 \hal^4 -240 \hal^3 + 660 \hal^2 - 720 \hal + 251 = - (20 \sqrt{6} -11)/6$, and
$- \hat{A}^2/2 + \hat{A} - k/2 + 1/12 = -(k-1)\{1 + (7+2\sqrt{6}) k\}/12$, where $\hat{A} := k \hal$.

The equality \eqref{ahat-suprisk} is directly obtained from \eqref{ahat-risk} because $1/(6 \sqrt{6}) - 11/720 > 0$.
\qed

\vspace{0.5cm}

We see that
the Bayesian predictive density $p_{\pi_\mathrm{J}}(y \mid x)$ based on the Jeffreys prior $\pi_{\mathrm J}$ is not asymptotically minimax.
The Jeffreys prior $\pi_\mathrm{J}$ is a Dirichlet prior $\bpial$ with $\alpha = 1/2$.
Thus, $\displaystyle 6 \alpha^2 - 12 \alpha + 5 = 1/2$, and
$- A^2/2 + A - k/2 + 1/12 = - (3 k^2 -2)/24$, where $A = k\alpha = k/2$.
Thus, from Theorem 1, we have
\begin{align*}
\sup_{\theta \in \Delta_{\epsN}} & \Biggl| R(\theta,p_{\pi_\mathrm{J}}(y \mid x)) - \frac{k-1}{2N}
- \frac{1}{N^2}
\biggl\{\sum_{i=1}^k \frac{1}{24 \theta_i}
-\frac{1}{24}(3 k^2 -2) \biggr\}
 \Biggr|
= \mathrm{o} (N^{-2} \epsN^{-1}).
\end{align*}
By putting $\theta_1 = \epsN$ and $\theta_i = (1-\epsN)/(k-1)$ $(i=2,\ldots,k)$, we have
\[
\sup_{\theta \in \Delta_{\epsN}} \left\{ R(\theta,p_{\pi_\mathrm{J}}(y \mid x)) - \frac{k-1}{2N} \right\}
\geq \frac{1}{24} \frac{1}{N^2 \epsN}
+ \mathrm{o} (N^{-2} \epsN^{-1}).
\]
Therefore, $p_{\pi_\mathrm{J}}(y \mid x)$ is not asymptotically minimax.

From Corollary 2, we obtain Corollary 3,
which is used to prove Theorem 3 in the next section.

We define
\begin{align*}
 \pi_a^{(N)} (\theta) =
\begin{cases}
\displaystyle \frac{\pi_a(\theta)}
{\displaystyle \int_{\Delta_\epsN} \pi_a(\theta)
\dd \theta_1 \cdots \dd \theta_{k-1}},&
 \theta \in \Delta_\epsN \\
0,& \text{otherwise},
\end{cases}
\end{align*}
and
\begin{align*}
 \bpi_\alpha^{(N)} (\theta) =
\begin{cases}
\displaystyle \frac{\bpi_\alpha(\theta)}
{\displaystyle \int_{\Delta_\epsN} \bpi_\alpha(\theta)
\dd \theta_1 \cdots \dd \theta_{k-1}},&
 \theta \in \Delta_\epsN \\
0,& \text{otherwise}.
\end{cases}
\end{align*}
The Bayes risk of a predictive density $q(y ; x)$ with respect to a prior $\pi$ is denoted by
\[
 R(\pi, q(y ; x)) := \int \pi(\theta) R(\theta, q(y ; x)) \dd \theta.
\]

\vspace{0.5cm}

\noindent
Corollary 3.
Suppose that $\{\epsN\}$ is a decreasing real number sequence such that
$\lim\limits_{N \rightarrow \infty} \epsN = 0$,
$\lim\limits_{N \rightarrow \infty} N^{\frac{3}{4}}\epsN = \infty$,
and $0 < \epsN < 1/k$ for every $N$.
Then,
\begin{align*}
R(\bpi_{\hal}^{(N)},p_{\bpihal}(y \mid x))
=& \sup_{\theta \in \Delta_{\epsN}} R(\theta,p_{\bpihal}(y \mid x))
+ \mathrm{o}(N^{-2}) \notag \\
=& \frac{k-1}{2N} - \frac{k-1}{12 N^2} \bigl\{1 + (7+2\sqrt{6})k \bigr\} + \mathrm{o}(N^{-2}).
\end{align*}
\qed

\vspace{0.5cm}

\noindent
Proof of Corollary 3.
From \eqref{ahat-risk}, we obtain
\begin{align*}
R(\bpi^{(N)}_{\hal}&,p_{\bpi_{\hal}}(y \mid x))=
\int_{\Delta_{\epsN}} \bpi_{\hal}^{(N)}(\theta) R(\theta,p_{\bpi_{\hal}}(y \mid x)) \dd \theta \\
=& \frac{k-1}{2N}
-\frac{k-1}{12 N^2}
\bigl\{1 + (7+2\sqrt{6})k \bigr\} \\
&+ \int_{\Delta_{\epsN}} \bpi_{\hal}^{(N)}(\theta) \biggl\{
- \frac{1}{N^3} \frac{1}{18 \sqrt{6}} \sum_{i=1}^k \frac{1}{\theta_i^2}
- \frac{1}{N^4} \left(\frac{1}{6\sqrt{6}} - \frac{11}{720} \right) \sum_{i=1}^k \frac{1}{\theta_i^3}
\biggr\} \dd \theta.
\end{align*}
Here, we have
\begin{align*}
\biggl| \int_{\Delta_{\epsN}} & \bpi_{\hal}^{(N)}(\theta) \biggl\{
- \frac{1}{N^3} \frac{1}{18 \sqrt{6}} \sum_{i=1}^k \frac{1}{\theta_i^2}
- \frac{1}{N^4} \left(\frac{1}{6\sqrt{6}} - \frac{11}{720} \right) \sum_{i=1}^k \frac{1}{\theta_i^3}
\biggr\} \dd \theta \biggr| \\
=& C
\int_{\Delta_{\epsN}} \bpi_{\hal}(\theta) \biggl\{
\frac{1}{N^3} \frac{1}{18 \sqrt{6}} \sum_{i=1}^k \frac{1}{\theta_i^2}
+ \frac{1}{N^4} \left(\frac{1}{6\sqrt{6}} - \frac{11}{720} \right) \sum_{i=1}^k \frac{1}{\theta_i^3}
\biggr\} \dd \theta \\
\leq& \frac{kC}{18 \sqrt{6} N^3} \int_{\theta_1 \geq \epsN} \bpi_{\hal}(\theta)
\frac{1}{\theta_1^2} \dd \theta_1
+ \frac{kC}{N^4} \left( \frac{1}{6\sqrt{6}} - \frac{11}{720} \right) \int_{\theta_1 \geq \epsN} \bpi_{\hal}(\theta) 
\frac{1}{\theta_1^3}
\dd \theta_1,
\end{align*}
where $C = 1/\int_{\Delta_\epsN} \bpi_{\hal}(\theta) \dd \theta_1 \cdots \dd \theta_{k-1}$.

Since the marginal density of $\theta_1$ of the Dirichlet prior $\bpi_{\hal}$ is
the Beta density
$\theta_1^{\hal-1} (1-\theta_1)^{(k-1)\hal-1}/B(\hal, (k-1) \hal)$,
we have
\begin{align*}
\biggl| &\int_{\Delta_{\epsN}} \bpi_{\hal}^{(N)}(\theta) \biggl\{
- \frac{1}{N^3} \frac{1}{18 \sqrt{6}} \sum_{i=1}^k \frac{1}{\theta_i^2}
- \frac{1}{N^4} \left(\frac{1}{6\sqrt{6}} - \frac{11}{720} \right) \sum_{i=1}^k \frac{1}{\theta_i^3}
\biggr\} \dd \theta \biggr| \\
&\leq \frac{kC}{N^3} \frac{1}{18 \sqrt{6}} \int_{\theta_1 \geq \epsN}
\frac{\theta_1^{\hal-1} (1-\theta_1)^{(k-1)\hal-1}}{B(\hal, (k-1)\hal)}
 \frac{1}{\theta_1^2} \dd \theta_1
 + \frac{kC}{N^4} \left( \frac{1}{6\sqrt{6}} - \frac{11}{720} \right) \int_{\theta_i \geq \epsN}
\frac{\theta_1^{\hal-1} (1-\theta_1)^{(k-1)\hal-1}}{B(\hal, (k-1)\hal)}
\frac{1}{\theta_1^3}
\dd \theta_1  \\
&\leq \frac{k C}{N^3} \frac{1}{18 \sqrt{6}} \frac{1}{2-\hal} \frac{1}{\epsN^{2-\hal}}
\frac{1}{B(\hal, (k-1)\hal)}
+ \frac{k C}{N^4} \left( \frac{1}{6\sqrt{6}} - \frac{11}{720} \right) \frac{1}{3-\hal} \frac{1}{\epsN^{3-\hal}}
\frac{1}{B(\hal, (k-1)\hal)} \\
&
= \mathrm{o}(N^{-2})
\end{align*}
because $\mathrm{O}(N^3 \epsN^{2-\hal})
\geq \mathrm{O}(N^2 (N \epsN)) \geq \mathrm{O}(N^2)$ and
$\mathrm{O}(N^4 \epsN^{3-\hal}) = \mathrm{O}((N^3 \epsN^{2-\hal}) (N \epsN)) \geq \mathrm{O}(N^2)$. \\
\qed

\vspace{0.5cm}

We use the following Lemmas 1--3 to prove Theorem 1.
The proofs of the lemmas are given in the appendix.

\vspace{0.5cm}

\noindent
Lemma 1.
For every nonnegative integer $m$ and every $x > -1$,
\[ 
\sum^{2m}_{i=1} (-1)^{i-1} \frac{1}{i} x^i + \frac{1}{2m+1} \frac{x^{2m+1}}{1+x}
\le
\log (1+x) \le \sum^{2m+1}_{i=1} (-1)^{i-1} \frac{1}{i} x^i.
\]
\qed

\vspace{0.5cm}

\noindent
Lemma 2.
Let $\mu_m(N,\theta)$ be the $m$-th central moments of the binomial distribution $\mathrm{Bi}(N,\theta)$
with index $N$ and parameter $\theta$.
Suppose that $\{\epsN\}$ be a decreasing sequence of real numbers such that
$\lim\limits_{N \rightarrow \infty} \epsN =0$,
$\lim\limits_{N \rightarrow \infty} N \epsN = \infty$, and $0 < \epsN < 1$ for every $N$.
\begin{description}
\item{(1)~}
For every positive integer $l$,
there exists a positive constant $C_{2l-1}$ such that
$\displaystyle \frac{|\mu_{2l-1}(N,\theta)|}{(N \theta)^{l-1}} \le C_{2l-1}$ for all $\theta \in [\epsN, 1]$ and $N$.
\item{(2)~}
For every positive integer $l$,
there exists a positive constant $C_{2l}$ such that
$\displaystyle \frac{|\mu_{2l}(N,\theta)|}{(N \theta)^l} \le C_{2l}$ for all $\theta \in [\epsN, 1]$ and $N$.
\end{description}
\qed

\vspace{0.5cm}

\noindent
Lemma 3. 
Let $x$ be a random variable distributed according to the binomial distribution $\mathrm{Bi}(N,\theta)$.
Define
\[
 w = \frac{x-N\theta}{N\theta + a},
\]
where $a$ is a positive real number.
Suppose that $\{\epsN\}$ be a decreasing sequence of real numbers such that
$\lim\limits_{N \rightarrow \infty} \epsN =0$,
$\lim\limits_{N \rightarrow \infty} N \epsN = \infty$, and $0 < \epsN < 1$ for every $N$.
Then, for every nonnegative integer $l$,
there exists a constant $C^{(w)}_{2l+1}$ such that
\begin{align*}
\E_\theta \left( -\frac{w^{2l+1}}{1+w} \right)
\le  \frac{1}{(N\theta)^l} C^{(w)}_{2l+1}
\end{align*}
for all $\theta \in [\epsN, 1]$ and $N$.
\qed

\vspace{0.5cm}

By using the lemmas, we prove Theorem 1.

\vspace{0.5cm}

\noindent
Proof of Theorem 1.

From \eqref{60-1} and Lemma 1, we have
\begin{align}
R(\theta, p_{\pi_a}(y \mid x))
& \le \sum^k_{i=1} \theta_i \mathrm{E}_\theta \left\{\sum_{l=1}^{10} \frac{1}{l} (- w_i)^l
- \frac{1}{11} \frac{{w_i}^{11}}{1 + w_i} \right\}
+ \sum^k_{i=1} \theta_i \left\{\sum_{l=1}^{4} \frac{1}{l} (- s_i)^l
- \frac{1}{5} \frac{{s_i}^{5}}{1 + s_i} \right\},
\label{LB0}
\end{align}
and
\begin{align}
R(\theta, p_{\pi_a}(y \mid x))
& \ge \sum^k_{i=1} \theta_i \mathrm{E}_\theta \left\{\sum_{l=1}^{9} \frac{1}{l} (- w_i)^l \right\}
+ \sum^k_{i=1} \theta_i \left\{\sum_{l=1}^{5} \frac{1}{l} (- s_i)^l \right\}.
\label{UB0}
\end{align}

From Lemma 2, we have
\begin{align}
|\E_{\theta} \left( w_i^{2l-1} \right)|
= \frac{|\mu_{2l-1}(N,\theta_i)|}{(N\theta_i + a_i)^{2l-1}}
& \leq \frac{|\mu_{2l-1}(N,\theta_i)|}{(N \theta_i)^{2l-1}} \leq \frac{C_{2l-1}}{(N \theta_i)^{l}},
\label{wineq2}
\end{align}
and
\begin{align}
|\E_{\theta} \left( w_i^{2l} \right)|
= \frac{|\mu_{2l}(N,\theta_i)|}{(N\theta_i + a_i)^{2l}}
& \leq \frac{|\mu_{2l}(N,\theta_i)|}{(N \theta_i)^{2l}} \leq \frac{C_{2l}}{(N \theta_i)^l},
\label{wineq1}
\end{align}
for every $a_i > 0$.

Obviously, the inequality
\begin{align}
 \frac{1}{1+s_i} =& \frac{N \theta_i + A \theta_i}{N \theta_i + a_i} \leq \frac{A}{a_i}
 \label{ineq-s}
\end{align}
holds since $0 < \theta_i < 1$ and $0 < a_i < A$.

From \eqref{LB0}, \eqref{wineq2}, \eqref{wineq1}, \eqref{ineq-s},
$\mathrm{E}_{\theta_i} (w_i)  = 0$, $\sum_{i=1}^k \theta_i s_i = 0$, and Lemma 3, we have
\begin{align}
R(\theta, p_{\pi_a}(y \mid x)) 
\le& \sum^k_{i=1} \theta_i \mathrm{E}_{\theta} \left\{\sum_{l=2}^{8} \frac{1}{l} (- w_i)^l
+ \sum_{l=9}^{10} \frac{1}{l} (- w_i)^l \right\}
+ \sum_{i=1}^k \frac{C_{11}^{(w)}}{N^5 {\theta_i}^4} \notag \\
& + \sum_{i=1}^k \theta_i \sum_{l=2}^4 \frac{1}{l}(-s_i)^l
+ \frac{1}{5} \sum_{i=1}^k \theta_i \frac{A}{a_i} \left| \frac{a_i - A \theta_i}{N \theta_i + A \theta_i} \right|^5 \notag \\
\le& \sum^k_{i=1} \theta_i \mathrm{E}_{\theta} \left\{\sum_{l=2}^{8} \frac{1}{l} (- w_i)^l \right\}
+ \sum_{i=1}^k \theta_i \sum_{l=2}^4 \frac{1}{l}(-s_i)^l
+ \frac{C'}{N^5 {\varepsilon_N}^4},
\label{UB-2}
\end{align}
where $C'$ is a positive constant not depending on $N$ or $\theta$.

In a similar way,
from \eqref{UB0} and \eqref{wineq2}, we have
\begin{align}
R(\theta, & \; p_{\pi_a}(y \mid x)) \ge
\sum^k_{i=1} \theta_i \mathrm{E}_{\theta} \left\{\sum_{l=2}^{8} \frac{1}{l} (- w_i)^l \right\}
+ \sum_{i=1}^k \theta_i \sum_{l=2}^4 \frac{1}{l}(-s_i)^l
- \frac{C''}{N^5 {\varepsilon_N}^4},
\label{LB-2}
\end{align}
where $C''$ is a positive constant not depending on $N$ or $\theta$.

The first to eighth central moments of the binomial distribution $\mathrm{Bi}(N,\theta)$ are given by
\begin{align}
 \mu_1(N,\theta) =& 0, ~~~ \mu_2(N,\theta) = N\theta(1-\theta), ~~ \mu_3(N,\theta) = N\theta (1-\theta) (1-2\theta), \notag \\
 \mu_4(N,\theta) =& 3 N^2 \theta^2 (1 - \theta)^2 + N \theta (1 - \theta) (1 - 6 \theta + 6\theta^2), \notag \\
 \mu_5(N,\theta) =& 10 N^2 \theta^2 (1-\theta)^2  (1-2\theta) +  N \theta (1-\theta) (1-2\theta) (1 - 12 \theta + 12 \theta^2), \notag \\
 \mu_6(N,\theta) =& 15 N^3 \theta^3 (1-\theta)^3 + 5 N^2 \theta^2 (1-\theta)^2 \left(5 - 26 \theta + 26 \theta^2 \right)
 + N \theta \phi_{6,1}(\theta), \notag \\
 \mu_7(N,\theta) =& 105 N^3 \theta^3 (1-\theta)^3 (1-2\theta)
+ N^2 \theta^2 \phi_{7,2}(\theta) + N \theta \phi_{7,1}(\theta), \notag \\
 \mu_8(N,\theta) =& 105 N^4 \theta^4 (1 - \theta)^4 + N^3 \theta^3 \phi_{8,3}(\theta)
+ N^2 \theta^2 \phi_{8,2}(\theta) + N \theta \phi_{8,1}(\theta),
\label{moments}
 \end{align}
where $\phi_{i,j} (\theta)$ $((i,j)=(6,1), (7,1), (7,2), (8,1), (8,2), (8,3))$ are polynomials of $\theta$.

Therefore, by using \eqref{UB-2}, \eqref{LB-2}, \eqref{moments},
and the inequalities
\begin{align*}
\frac{1}{N\theta_i} &
\sum_{l=1}^{2m-1}
\left(- \frac{a_i}{N \theta_i}\right)^{l} 
=
\frac{1}{N \theta_i} \frac{\displaystyle 1 - \left(\frac{a_i}{N \theta_i}\right)^{2m}}{\displaystyle \displaystyle 1 + \frac{a_i}{N \theta_i}}
\le
 \frac{1}{N \theta_i + a_i}
\le
\frac{1}{N\theta_i}
\sum_{l=1}^{2m}
\left(- \frac{a_i}{N \theta_i}\right)^{l}
=
\frac{1}{N \theta_i} \frac{\displaystyle 1 + \left(\frac{a_i}{N \theta_i}\right)^{2m+1}}{\displaystyle \displaystyle 1 + \frac{a_i}{N \theta_i}},
\end{align*}
we obtain \eqref{theorem1} by a straightforward but lengthy calculation.
In addition to the calculation by hand, the result is verified by using a computer algebra software.
\qed
\section{Minimax predictive densities}

In this section, we prove that the Bayesian predictive density based on a Dirichlet prior $\bpi_{\hal}$,
where $\hal := 1+1/\sqrt{6}$,
is asymptotically minimax in the sense of \eqref{1-3-1} if $\{\epsN\}$ satisfies appropriate conditions.

The Bayesian predictive density 
with respect to the prior $\bpi_\alpha$ is given by
\begin{align}
p_{\bpi_\alpha} (y \mid x) = \frac{B(x_1+y_1+\alpha, \ldots, x_k+y_k+\alpha)}{B(x_1+\alpha, \ldots, x_k+\alpha)}
= \frac{\sum_{i=1}^{k} x_i y_i + \alpha}{N + k \alpha},
\label{ourpd}
\end{align}
and that with respect to the prior $\bpi_\alpha^{(N)}$ is given by
\begin{align}
p_{\bpi^{(N)}_\alpha} (y|x)
=& \frac{B_{\Delta_{\epsN}}(x_1+y_1+\alpha, \ldots, x_k+y_k+\alpha)}{B_{\Delta_{\epsN}}(x_1+\alpha, \ldots, x_k+\alpha)} \notag \\
=& \frac{\sum_{i=1}^{k} x_i y_i + \alpha}{N + k \alpha}
\frac{I_{\Delta_{\varepsilon_N}}(x_1+y_1+\alpha, \ldots, x_k+y_k+\alpha)}{I_{\Delta_{\varepsilon_N}}(x_1+\alpha, \ldots, x_k+\alpha)},
\label{BIdecomp}
\end{align}
where we define
\[
B_{\Delta_\varepsilon}(\alpha_1,\ldots,\alpha_k)
:= \int_{\Delta_\varepsilon} \theta_1^{\alpha_1-1} \cdots \theta_k^{\alpha_k-1}
\dd \theta_1 \cdots \dd \theta_{k-1}
\]
and
\[
 I_{\Delta_\varepsilon}(\alpha_1,\ldots,\alpha_k)
:= \frac{B_{\Delta_\varepsilon}(\alpha_1,\ldots,\alpha_k)} 
{B(\alpha_1,\ldots,\alpha_k)}
\]
for $\alpha_i > 0$ $(i=1,\ldots,k)$ and $0 < \varepsilon <1/k$.
If $k=2$, $I_{\Delta_\varepsilon}(\alpha_1, \alpha_2)
= \{\int_\eps^{1-\eps} \theta^{\alpha_1-1}(1-\theta)^{\alpha_2-1} \dd \theta\}
/ B(\alpha_1,\alpha_2)$.

In the proof of minimaxity of prediction,
the inequalities
\begin{align}
\sup_{\theta \in \Delta_{\epsN}} & R(\theta, p_{\pi}(y \mid x))
\geq \inf_{q} \sup_{\theta \in \Delta_{\epsN}} R(\theta, q(y ; x))
= \inf_{q} \sup_{\pi' \in \mathcal{P}(\Delta_{\epsN})} R(\pi', q(y ; x)) \notag \\
\geq& \sup_{\pi' \in \mathcal{P}(\Delta_{\epsN})} \inf_{q} R(\pi', q(y ; x))
\geq \inf_{q} R(\pi^*, q(y ; x))
= R(\pi^*, p_{\pi^*}(y \mid x)),
\label{fundamental-ineq}
\end{align}
which hold for every $\pi \in \mathcal{P}(\Delta)$ and $\pi^* \in \mathcal{P}(\Delta_{\epsN})$, play an essential role;
see \citet{GD04} for related inequalities in a very general setting.
Each inequality in \eqref{fundamental-ineq} is easy to verify.
The last inequality in \eqref{fundamental-ineq} is due to the fact, proved by \citet{Aitchison75},
that the Bayes risk of a predictive density with respect to a prior $\pi^*$ is minimized when it is
the Bayesian predictive density $p_{\pi^*}(y \mid x)$ based on $\pi^*$.
Thus, by putting $\pi^* = \bpi^{(N)}_\alpha$ in \eqref{fundamental-ineq}, we have
\[
R(\bpi^{(N)}_\alpha, p_{\bpi_\alpha}(y \mid x))
\ge \inf_{q} R(\bpi^{(N)}_\alpha, q(y ; x))
=  R(\bpi^{(N)}_\alpha, p_{\bpi^{(N)}_\alpha}(y \mid x)).
\]

In the following, we first prove Theorem 2 that shows that the difference
$R(\bpi^{(N)}_\alpha, p_{\bpi_\alpha}(y \mid x)) - R(\bpi^{(N)}_\alpha, p_{\bpi^{(N)}_\alpha}(y \mid x))$
is $\mathrm{O}(N^{-1}\epsN^{\alpha})$ if $\{\epsN\}$ satisfies appropriate conditions.
Next, combining Corollary 3 and Theorem 2, we prove Theorem 3 showing that $p_{\pi^{(N)}_a}(y \mid x)$
is asymptotically minimax under suitable conditions.

\vspace{0.5cm}

\noindent
Theorem 2.
Let $p_{\bpi^{(N)}_\alpha}(y \mid x)$ and $p_{\bpi_\alpha}(y \mid x)$ be predictive densities \eqref{ourpd} and \eqref{BIdecomp},
respectively.
Suppose that $\{ \epsN \}$ is a decreasing sequence of real numbers such that
$\lim\limits_{N \rightarrow \infty} \epsN = 0$,
$\lim\limits_{N \rightarrow \infty} N \epsN = \infty$,
and $0 < \epsN < 1/k$ for every $N$.
Then the difference of the Bayes risks of $p_{\bpi^{(N)}_\alpha}(y \mid x)$ and $p_{\bpi_\alpha}(y \mid x)$ with respect to $\bpi_\alpha^{(N)}$
satisfies
\begin{align*}
R(\bpi^{(N)}_\alpha, p_{\bpi_\alpha}(y \mid x))
- R(\bpi^{(N)}_\alpha, p_{\bpi^{(N)}_\alpha}(y \mid x))
= \mathrm{O}(N^{-1} \epsN^{\alpha}).
\end{align*}
\qed

\vspace{0.5cm}

Theorem 2 means that the disadvantage of adopting a prior $\bpi_\alpha$ that
does not satisfy \\ $\int_{\Delta_\epsN} \bpi_\alpha(\theta) \dd \theta_1 \cdots \dd \theta_{k-1} = 1$
is asymptotically small.

We use Lemmas 4--8 below to prove Theorem 2.
The proofs of the lemmas are given in the Appendix.

\vspace{0.5cm}

\noindent
Lemma 4.
For every $\alpha_1 > 0$, \ldots, $\alpha_k > 0$ and $0 < \varepsilon < 1/k$,
\begin{align*}
I_{\Delta_\varepsilon} (\alpha_1 + 1, \alpha_2, \dotsb, \alpha_k) - I_{\Delta_\varepsilon} (\alpha_1, \dotsb, \alpha_k)
\leq& \;
\frac{\Gamma (\sum_{i=1}^{k} \alpha_i)}{\Gamma(\alpha_1 + 1) \Gamma (\sum_{i=2}^{k} \alpha_i)}
\eps^{\alpha_1} (1-\eps)^{\alpha_2+\dotsb+\alpha_k}.
\end{align*}

\vspace{0.5cm}

\noindent
Lemma 5.
If $0 \leq s < t \leq 1$, $0 \leq u < v \leq 1$, $s \leq u$, and $t \leq v$, then for all $\alpha>0$ and $\beta>0$,
\begin{align*}
\frac{B_{[s,t]}(\alpha+1,\beta)}{B_{[s,t]}(\alpha,\beta)}
\leq \frac{B_{[u,v]}(\alpha+1,\beta)}{B_{[u,v]}(\alpha,\beta)},
\end{align*}
where
\[
 B_{[s,t]}(\alpha,\beta) := \int_s^t \theta^{\alpha-1} (1-\theta)^{\beta-1} \dd \theta.
\]
\qed

\vspace{0.5cm}

\noindent
Lemma 6.
For every $\alpha_1 > 0, \ldots, \alpha_k > 0$, and $0 < \varepsilon < 1/k$, the inequality
\begin{align*}
\frac{B_{\Delta_\varepsilon}(\alpha_1+1,\alpha_2,\ldots,\alpha_k)}
{B_{\Delta_\varepsilon}(\alpha_1,\alpha_2,\ldots,\alpha_k)} \leq
\frac{B_{[\varepsilon,1]}(\alpha_1+1,\alpha_2+\cdots+\alpha_k)}{B_{[\varepsilon,1]}(\alpha_1,\alpha_2+\cdots+\alpha_k)}
\end{align*}
holds.
\qed

\vspace{0.5cm}

\noindent
Lemma 7. 
For every $\alpha_1 > 0, \ldots, \alpha_k > 0$,
the equality
\begin{align*}
\sum_{x_2,\ldots,x_{k}: \sum_{i=2}^{k} x_i = N-x_1} &
\frac{B(x_1+\alpha_1,x_2+\alpha_2,\ldots,x_k+\alpha_k)}{B(\alpha_1,\alpha_2,\ldots,\alpha_k)}{\binom{N}{x_1, \ldots, x_k}} \\
=&
\frac{B(x_1+\alpha_1,N-x_1 + \sum_{i=2}^k \alpha_i)}{B(\alpha_1, \sum_{i=2}^k \alpha_i)}{\binom{N}{x_1}}
\end{align*}
holds.
\qed

\vspace{0.5cm}

\noindent
Lemma 8. 
For every $\alpha >0$, $\beta >0$, and $\varepsilon \in [0,1)$, the inequality
\begin{align*}
& \frac{B_{[\varepsilon,1]} (\alpha+1, \beta)}
{B_{[\varepsilon,1]} (\alpha, \beta)}
= \frac{\int^1_\varepsilon  \theta^\alpha (1-\theta)^{\beta-1} \dd \theta}
{\int^1_\varepsilon \theta^{\alpha-1} (1-\theta)^{\beta-1} \dd \theta}
= \frac{\alpha}{\alpha+\beta} + \frac{1}{\alpha+\beta}
\frac{\varepsilon^\alpha (1-\varepsilon)^\beta}
{B_{[\varepsilon,1]} (\alpha, \beta)}
\le \frac{(1-\eps)\alpha}{\alpha+\beta} + \eps
\end{align*}
holds.
\qed

\vspace{0.5cm}

By using the lemmas, we prove Theorem 2.

\vspace{0.5cm}

\noindent
Proof of Theorem 2.
From \eqref{ourpd} and \eqref{BIdecomp}, the difference between the risk functions of
$R(\theta, p_{\bpi_\alpha} (y \mid x))$ and $R(\theta, p_{\bpi_\alpha^{(N)}} (y \mid x))$ is given by
\begin{align*}
R(&\theta, p_{\bpi_\alpha} (y \mid x)) - R(\theta, p_{\bpi_\alpha^{(N)}} (y \mid x))
= \sum_x \sum_y p(x,y|\theta) \log \frac{p_{\bpi^{(N)}_\alpha} (y|x)}{p_{\bpi_\alpha} (y|x)} \notag \\
=& \sum_x \sum_y \theta_1^{y_1} \cdots \theta_k^{y_k}{\binom{N}{x_1, \ldots, x_k}} \theta_1^{x_1} \theta_2^{x_2} \cdots \theta_k^{x_k}
\log \frac{I_{\Delta_{\epsN}}(x_1+y_1+\alpha, x_2+y_2+\alpha, \ldots, x_k+y_k+\alpha)}{I_{\Delta_\epsN}(x_1+\alpha, \ldots, x_k+\alpha)}.
\end{align*}

To evaluate the difference between the Bayes risks $R(\bpi^{(N)}_\alpha, p_{\bpi_\alpha} (y \mid x))$
and $R(\bpi^{(N)}_\alpha, p_{\bpi_\alpha^{(N)}} (y \mid x))$,
it is sufficient to consider the case $y_1 = 1$
because of the symmetry of the index $i$.
Thus,
\begin{align*}
R(&\bpi^{(N)}_\alpha, p_{\bpi_\alpha} (y \mid x)) - R(\bpi^{(N)}_\alpha, p_{\bpi_\alpha^{(N)}} (y \mid x)) \notag \\
=& k \int_{\Delta_\epsN}
\frac{{\theta_1}^{\alpha-1} \cdots {\theta_k}^{\alpha-1}}{B_{\Delta_\epsN} (\alpha,\alpha,\ldots,\alpha)}
\theta_1 \sum_x {\binom{N}{x_1, \ldots, x_k}} \theta_1^{x_1} \theta_2^{x_2} \cdots \theta_k^{x_k} \notag \\
& \times \log \frac{I_{\Delta_\epsN}(x_1+1+\alpha, x_2+\alpha, \ldots, x_k+\alpha)}
{I_{\Delta_\epsN}(x_1+\alpha, \ldots, x_k+\alpha)} \dd \theta_1 \cdots \dd \theta_{k-1} \notag \\
=& \frac{k}{B_{\Delta_\epsN} (\alpha,\alpha,\ldots,\alpha)}
\sum_x B_{\Delta_\epsN}(x_1+1+\alpha,x_2+\alpha,\ldots,x_k+\alpha)
\binom{N}{x_1, \ldots, x_k} \notag \\
& \times \log \frac{I_{\Delta_\epsN}(x_1+1+\alpha, x_2+\alpha, \ldots, x_k+\alpha)}{I_{\Delta_\epsN}(x_1+\alpha, \ldots, x_k+\alpha)}.
\end{align*}
Because $\log(x+1) \leq x$ for $x > -1$, we have
\begin{align*}
R(&\bpi^{(N)}_\alpha, p_{\bpi_\alpha} (y \mid x)) - R(\bpi^{(N)}_\alpha, p_{\bpi_\alpha^{(N)}} (y \mid x)) \notag \\
\leq& \frac{k}{B_{\Delta_\epsN} (\alpha,\alpha,\ldots,\alpha)}
\sum_x B_{\Delta_\epsN}(x_1+1+\alpha,x_2+\alpha,\ldots,x_k+\alpha)
\binom{N}{x_1, \ldots, x_k} \notag \\
& \times \left\{ \frac{I_{\Delta_\epsN}(x_1+1+\alpha, x_2+\alpha, \ldots, x_k+\alpha)}{I_{\Delta_\epsN}(x_1+\alpha, \ldots, x_k+\alpha)}
- 1 \right\}
\notag \\
=& \frac{k}{B_{\Delta_\epsN} (\alpha,\alpha,\ldots,\alpha)}
\sum_x B(x_1+\alpha,\ldots,x_k+\alpha) \notag \\
& \times \frac{B_{\Delta_\epsN}(x_1+1+\alpha, x_2+\alpha, \ldots, x_k+\alpha)}{B_{\Delta_\epsN}(x_1+\alpha, \ldots, x_k+\alpha)}
\binom{N}{x_1, \ldots, x_k} \notag \\
& \times \bigl\{ I_{\Delta_\epsN}(x_1+1+\alpha, x_2+\alpha, \ldots, x_k+\alpha) - I_{\Delta_\epsN}(x_1+\alpha, \ldots, x_k+\alpha) \bigr\}.
\end{align*}
From Lemma 4, we obtain
\begin{align*}
R(&\bpi^{(N)}_\alpha, p_{\bpi_\alpha} (y \mid x)) - R(\bpi^{(N)}_\alpha, p_{\bpi_\alpha^{(N)}} (y \mid x)) \notag \\
\leq& \frac{k}{B_{\Delta_\epsN} (\alpha,\alpha,\ldots,\alpha)}
\sum_x B(x_1+\alpha,\ldots,x_k+\alpha)
\frac{B_{\Delta_\epsN}(x_1+1+\alpha, x_2+\alpha, \ldots, x_k+\alpha)}{B_{\Delta_\epsN}(x_1+\alpha, \ldots, x_k+\alpha)} \notag \\
& \times \binom{N}{x_1, \ldots, x_k}
\frac{\Gamma(N+k\alpha)}{\Gamma(x_1+1+\alpha) \Gamma(N-x_1+(k-1)\alpha)}
\epsN ^{x_1+\alpha} (1-\epsN)^{N-x_1+(k-1)\alpha}.
\end{align*}
From Lemmas 6 and 7, we have
\begin{align*}
R(&\bpi^{(N)}_\alpha, p_{\bpi_\alpha} (y \mid x)) - R(\bpi^{(N)}_\alpha, p_{\bpi_\alpha^{(N)}} (y \mid x)) \notag \\
\leq& \frac{k}{B_{\Delta_\epsN} (\alpha,\alpha,\ldots,\alpha)}
\sum_x B(x_1+\alpha,\ldots,x_k+\alpha)
\frac{B_{[\epsN,1]}(x_1+1+\alpha, N - x_1 + (k-1)\alpha)}
{B_{[\epsN,1]}(x_1+\alpha, N - x_1 + (k-1)\alpha)} \notag \\
& \times \binom{N}{x_1, \ldots, x_k}
\frac{\Gamma(N+k\alpha)}{\Gamma(x_1+1+\alpha) \Gamma(N-x_1+(k-1)\alpha)} \epsN^{x_1+\alpha} (1-\epsN)^{N-x_1+(k-1)\alpha}
\notag \\
=& \frac{k B(\alpha,\alpha,\ldots,\alpha)}{B_{\Delta_\epsN} (\alpha,\alpha,\ldots,\alpha) B(\alpha,(k-1)\alpha)}
\sum_{x_1=0}^N \frac{1}{x_1+\alpha}
\frac{B_{[\epsN,1]}(x_1+1+\alpha, N-x_1+(k-1)\alpha)}
{B_{[\epsN,1]}(x_1+\alpha, N-x_1+(k-1)\alpha)} \notag \\
& \times \binom{N}{x_1}
\epsN^{x_1+\alpha} (1-\epsN)^{N-x_1+(k-1)\alpha}.
\end{align*}
Since
\[
 \frac{B_{[\epsN,1]}(x_1+1+\alpha, N-x_1+(k-1)\alpha)}
{B_{[\epsN,1]}(x_1+\alpha, N-x_1+(k-1)\alpha)}
\leq \frac{(1-\epsN)(x_1+\alpha)}{N+k\alpha} + \epsN
\]
because of Lemma 8,
we have
\begin{align*}
R(&\bpi^{(N)}_\alpha, p_{\bpi_\alpha} (y \mid x)) - R(\bpi^{(N)}_\alpha, p_{\bpi_\alpha^{(N)}} (y \mid x)) \notag \\
\le& \frac{k \epsN^{\alpha} (1-\epsN)^{(k-1)\alpha}}{I_{\Delta_\epsN} (\alpha,\ldots,\alpha) B(\alpha,(k-1)\alpha)}
\left\{ \frac{1-\epsN}{N+k\alpha} 
+ \epsN \sum_{x_1=0}^N  \frac{1}{x_1+\alpha}
\binom{N}{x_1} \epsN^{x_1} (1-\epsN)^{N-x_1} \right\} \notag \\
=& \frac{k \epsN^{\alpha} (1-\epsN)^{(k-1)\alpha}}{I_{\Delta_\epsN} (\alpha,\ldots,\alpha) B(\alpha,(k-1)\alpha)}
\left\{ \frac{1-\epsN}{N+k\alpha} \right. \notag \\
&\left. + \sum_{x_1=0}^N \frac{x_1+1}{x_1+\alpha} \frac{1}{N+1}
\frac{(N+1)!}{(x_1+1)! (N+1-x_1-1)!} \epsN^{x_1+1} (1-\epsN)^{N+1-x_1-1} \right\} \notag \\
=& \frac{k \epsN^{\alpha} (1-\epsN)^{(k-1)\alpha}}{I_{\Delta_\epsN} (\alpha,\ldots,\alpha) B(\alpha,(k-1)\alpha)}
\left\{ \frac{1-\epsN}{N+k\alpha} \right.
\left. + \sum_{z=0}^{N+1} \frac{z}{z+\alpha-1} \frac{1}{N+1}
\binom{N+1}{z} \epsN^{z} (1-\epsN)^{N+1-z} \right\}, \notag
\end{align*}
where we define $z/(z+\alpha-1) = 0$ if $\alpha=1$ and $z=0$.

Since there exists a constant $C_\alpha > 0$ such that
$|z/(z+\alpha-1)| < C_\alpha$ for every $z$,
we have
\begin{align*}
R(\bpi^{(N)}_\alpha,& \; p_{\bpi_\alpha} (y \mid x)) - R(\bpi^{(N)}_\alpha, p_{\bpi_\alpha^{(N)}} (y \mid x)) \notag \\
\leq& \, \frac{k \epsN^{\alpha} (1-\epsN)^{(k-1)\alpha}}{I_{\Delta_\epsN} (\alpha,\ldots,\alpha) B(\alpha,(k-1)\alpha)}
\left( \frac{1-\epsN}{N+k\alpha}
+ \frac{C_\alpha}{N+1} \right)
= \mathrm{O}(N^{-1} \epsN^{\alpha}).
\end{align*}
\qed

\vspace{0.5cm}

Now we prove Theorem 3 that shows $p_{\bpi_{\hal}} (y \mid x)$,
where $\hal = 1 + 1/\sqrt{6}$, is asymptotically minimax.
The constant $1/\hal = \sqrt{6}/(\sqrt{6}+1)$ in the theorem is approximately $0.7101$.

\vspace{0.5cm}

\noindent
Theorem 3.
Let $p_{\bpi_{\hal}}(y \mid x)$ be the predictive density based on the prior
\[
 \bpi_{\hal}(\theta) \dd \theta_1 \cdots \dd \theta_{k-1}
\propto
{\theta_1}^{1/\sqrt{6}} \cdots {\theta_{k-1}}^{1/\sqrt{6}}
(1-\sum_{i=1}^{k-1} \theta_i)^{1/\sqrt{6}}
 \dd \theta_1 \cdots \dd \theta_{k-1}.
\]
Suppose that $\{\epsN\}$ be a decreasing sequence of real numbers such that
$\lim\limits_{N \rightarrow \infty} N^{3/4} \epsN = \infty$,
$\lim\limits_{N \rightarrow \infty} N^{1/\hal} \epsN = 0$,
and $0< \epsN < 1/k$ for every $N$.
Then, 
\[
\sup_{\theta \in \Delta_{\epsN}} R(\theta, p_{\bpi_{\hal}}(y \mid x))
= \inf_{q} \sup_{\theta \in \Delta_{\epsN}} R(\theta, q(y ; x))
+ \mathrm{o}(N^{-2}).
\]
\qed

\vspace{0.5cm}

\noindent
Proof.
By setting $\pi = \bpi_{\hat{\alpha}}$ and $\pi^* = \bpi^{(N)}_{\hat{\alpha}}$ in \eqref{fundamental-ineq},
we obtain
\begin{align}
\sup_{\theta \in \Delta_{\epsN}} R(\theta, p_{\bpi_{\hat{\alpha}}}(y \mid x))
\geq& \inf_q \sup_{\theta \in \Delta_\epsN} R(\theta, q(y ; x))
\geq R(\bpi^{(N)}_{\hat{\alpha}}, p_{\bpi_{\hat{\alpha}}^{(N)}}(y \mid x)).
\label{fundamental-ineq2}
\end{align}
From Theorem 2, we have
\begin{align}
R(\bpi^{(N)}_{\hat{\alpha}}, p_{\bpi_{\hat{\alpha}}}(y \mid x))
- R(\bpi^{(N)}_{\hat{\alpha}}, p_{\bpi^{(N)}_\alpha}(y \mid x))
= \mathrm{O}(N^{-1} \epsN^{\hal}) = \mathrm{o}(N^{-2}),
\label{fundamental-ineq3}
\end{align}
because $\epsN = \mathrm{o}(N^{-1/\hal})$.
From \eqref{fundamental-ineq2} and \eqref{fundamental-ineq3}, we have
\begin{align}
\sup_{\theta \in \Delta_{\epsN}} R(\theta, p_{\bpi_{\hat{\alpha}}}(y \mid x))
\geq& \inf_q \sup_{\theta \in \Delta_\epsN} R(\theta, q(y ; x))
\geq R(\bpi^{(N)}_{\hat{\alpha}}, p_{\bpi_{\hat{\alpha}}}(y \mid x))
+ \mathrm{o}(N^{-2}).
\label{fundamental-ineq4}
\end{align}
Here, from Corollary 3,
\begin{align}
R(\bpi_{\hal}^{(N)},p_{\bpihal}(y \mid x))
=& \sup_{\theta \in \Delta_{\epsN}} R(\theta,p_{\bpihal}(y \mid x))
+ \mathrm{o}(N^{-2}).
\label{however}
\end{align}
From \eqref{fundamental-ineq4} and \eqref{however}, we obtain the desired equality.
\qed
\section{Discussion}

The results in the present paper indicate that $\bpi_{\hat{\alpha}}(\theta) \propto
{\theta_1}^{1/\sqrt{6}} \cdots {\theta_{k-1}}^{1/\sqrt{6}}
(1-\sum_{i=1}^{k-1} \theta_i)^{1/\sqrt{6}}$
could be a reasonable objective prior for one-step ahead prediction.
The prior $\bpi_{\hal}(\theta)$ can be regarded as an asymptotic approximation to the latent information prior,
based on which a minimax predictive density is constructed, and it seems to consistent with some numerical results in \citet{Komaki11}.
Bayesian predictive densities based on commonly used objective priors, such as the Jeffreys priors $\pi_\mathrm{J}$ on $\Delta$,
$\pi_\mathrm{J}^{(N)}$ on $\Delta_\epsN$, the uniform priors $\pi_{\mathrm{U}}$ on $\Delta$,
or $\pi_{\mathrm{U}}^{(N)}$ on $\Delta_\epsN$ are not asymptotically minimax.

The conditions $\lim\limits_{N \rightarrow \infty} N^{3/4} \epsN = \infty$ and
$\lim\limits_{N \rightarrow \infty} N^{1/\hal} \epsN = 0$ assumed in Theorem 3
are sufficient conditions.
If $\epsN$ converges to $0$ very rapidly, then the condition $\lim\limits_{N \rightarrow \infty} N^{3/4} \epsN = \infty$
is not satisfied and we need to take into consideration the singularity at the boundary of the parameter space $\Delta$.
If $\epsN$ converges to $0$ very slowly, then the condition
$\lim\limits_{N \rightarrow \infty} N^{1/\hal} \epsN = 0$
is not satisfied and we cannot neglect the difference between $\bpi_{\hat{\alpha}}$ and $\bpi^{(N)}_\alpha$.
The constant $1/\hal = \sqrt{6}/(\sqrt{6}+1) \simeq 0.7101$ is not much smaller than 3/4.
It may be possible to weaken the condition
$\lim\limits_{N \rightarrow \infty} N^{3/4} \epsN = \infty$ by using an expansion of the risk function
with more higher order terms than the formula
in Theorem 1.

\newpage
\begin{appendix}
\section{Proofs of lemmas}

\noindent
Proof of Lemma 1.
(1)
Let
\[
f(x) := \log (1+x) - \sum^{2m+1}_{i=1} (-1)^{i-1} \frac{1}{i} x^i.
\]
Then, $f(0) = 0$, and
\begin{align*}
f'(x) =& \frac{1}{1+x} - \sum^{2m}_{i=0} (-1)^{i} x^{i}
= \frac{1}{1+x} - \left(\frac{1}{1+x} + \frac{x^{2m+1}}{1+x}\right) = - \frac{x^{2m+1}}{1+x}.
\end{align*}
Thus, $f'(x) > 0$ for $-1 < x < 0$, $f'(x) = 0$ for $x = 0$, and $f'(x) < 0$ for $x > 0$.
Therefore, $f(x) \le 0$ for $x > -1$, and the equality holds only when $x = 0$.

\noindent
(2)
Let
\[
f(x) := \log (1+x) - \sum^{2m}_{i=1} (-1)^{i-1} \frac{1}{i} x^i - \frac{1}{2m+1} \frac{x^{2m+1}}{1+x}.
\]
Then, $f(0) = 0$, and
\begin{align*}
f'(x) =& \frac{1}{1+x} - \sum^{2m-1}_{i=0} (-1)^{i} x^{i} - \frac{x^{2m}}{1+x} + \frac{1}{2m+1} \frac{x^{2m+1}}{(1+x)^2} \\
=& \frac{1}{1+x} - \left(\frac{1}{1+x} - \frac{x^{2m}}{1+x} \right)
- \frac{x^{2m}}{1+x} + \frac{1}{2m+1} \frac{x^{2m+1}}{(1+x)^2}
= \frac{1}{2m+1} \frac{x^{2m+1}}{(1+x)^2}
\end{align*}
Thus, $f'(x) < 0$ for $-1 < x < 0$, $f'(x) = 0$ for $x = 0$, and $f'(x) > 0$ for $x > 0$.
Therefore, $f(x) \ge 0$ for $-1 < x$, and the equality holds only when $x = 0$.
\qed

\vspace{0.5cm}

\noindent
Proof of Lemma 2.
We prove the desired results by induction.
Assume that $\mu_{2l-1} (N, \theta)$ and $\mu_{2l} (N, \theta)$,
where $l$ is a positive integer, are represented as
\begin{align}
\mu_{2l-1}(N,\theta) = \sum_{i=1}^{l-1} f_{2l-1,i}(\theta) (N \theta)^i \mbox{~~and~~}
\mu_{2l}(N,\theta) = \sum_{i=1}^l f_{2l,i}(\theta) (N \theta)^i
\label{recursive}
\end{align}
where $f_{2l-1,i}(\theta)$ $(i=1,2,\ldots,l-1)$ and $f_{2l,i}(\theta)$ $(i=1,2,\ldots,l)$
are polynomials with integer coefficients.
Then,
by using the recurrence equation
\begin{align*}
\mu_{m+1}(N,\theta)  =& \theta(1-\theta) \left\{ N m \mu_{m-1}(N,\theta) +\frac{\dd \mu_m (N,\theta)}{\dd \theta} \right\}
~~~~~ (m=2,3,4,\ldots)
\end{align*}
by \citet{Romanovsky23},
we have
\begin{align*}
\mu_{2l+1}(N,\theta) =& \theta (1-\theta) \left \{2 N l \sum_{i=1}^{l-1} f_{2l-1,i}(\theta) (N \theta)^i
+\frac{\dd }{\dd  \theta} \sum_{i=1}^l f_{2l,i}(\theta) (N \theta)^i \right\}
\end{align*}
and
\begin{align*}
\mu_{2l+2}(N,\theta) =& \theta(1-\theta) \left \{2 N (l+1) \sum_{i=1}^l f_{2l,i}(\theta) (N \theta)^i
+\frac{\dd }{\dd  \theta} \sum_{i=1}^l f_{2l+1,i}(\theta) (N \theta)^i \right\}.
\end{align*}
Thus, $\mu_{2l+1}(N,\theta)$ and $\mu_{2l+2}(N,\theta)$ are represented as
\[
\mu_{2l+1}(N,\theta)
= \sum_{i=1}^l f_{2l+1,i}(\theta) (N \theta)^i \mbox{~~and~~}
\mu_{2l+2}(N,\theta) = \sum_{i=1}^{l+1} f_{2l+2,i}(\theta) (N \theta)^i,
\]
where $f_{2l+1,i}(\theta)$ and $f_{2l+2,i}(\theta)$ are polynomials of $\theta$ with integer coefficients.

Since $\mu_1(N,\theta) = 0$ and $\mu_2(N,\theta) = N\theta(1-\theta)$,
the equation \eqref{recursive} holds for every positive integer $l$.

Therefore, because $N \epsN$ goes to infinity,
there exist constants $C_{2l-1}$ and $C_{2l}$ not depending on $N$ or $\theta$ such that
\begin{align*}
\frac{|\mu_{2l-1}(N,\theta)|}{(N \theta)^{l-1}}
\le& \sum_{i=1}^{l-1} \frac{|f_{2l-1,i}(\theta)|}{(N \theta)^{l-1-i}}
\le \sum_{i=1}^{l-1} \frac{\max_{\theta \in [0,1]} |f_{2l-1,i}(\theta)|}{(N \epsN)^{l-1-i}} \le C_{2l-1} \\
\intertext{and}
\frac{|\mu_{2l}(N,\theta)|}{(N \theta)^l} \le& \sum_{i=1}^l \frac{|f_{2l,i}(\theta)|}{(N \theta)^{l-i}}
\le \sum_{i=1}^l \frac{\max_{\theta \in [0,1]} |f_{2l,i}(\theta)|}{(N \epsN)^{l-i}} \le C_{2l},
\end{align*}
respectively.
\qed

\vspace{0.5cm}

\noindent
Proof of Lemma 3.
We have
\begin{align*}
 \mathrm{E} \biggl( & \frac{w^{2l}} { 1 + w} \biggr)
=  \sum^N_{x=0} \binom{N}{x} \theta^x (1 - \theta)^{N-x} \frac{N \theta + a}{x+a} \left( \frac{x - N \theta}{N \theta + a} \right)^{2l} \\
=& \frac{1}{(N \theta +a)^{2l-1}}
\sum^N_{x=0} \frac{(N+1)!}{ (x+1) ! \{ (N+1) - (x+1) \} !} 
\frac{1}{N+1} \frac{1}{\theta} \theta^{x+1} (1-\theta)^{N+1 - (x+1)}
\frac{x+1}{x+a}
(x - N \theta)^{2l}.
\end{align*}
Here, for every $x \geq 0$,
\[
\frac{x+1}{x+a}
\leq
\left\{
\begin{array}{cc}
1, & \mbox{if } a \geq 1, \\[5pt]
\displaystyle \frac{1}{a}, & \mbox{if } 0< a < 1.
\end{array}
\right.
\]
Thus,
\begin{align*}
 \mathrm{E} \biggl( & \frac{w^{2l}} {1 + w} \biggr)
\leq
\frac{\displaystyle \max \left( 1, \; \frac{1}{a} \right)}{(N \theta +a)^{2l-1} (N+1)\theta}
\sum^{N+1}_{z=0} \frac{(N+1)!}{ z! (N+1-z)!} \theta^z (1-\theta) ^{N+1-z} 
\{z - (N+1) \theta - (1-\theta)\}^{2l} \\
&= \frac{\displaystyle \max \left( 1, \; \frac{1}{a} \right)}{(N \theta +a)^{2l-1} (N+1)\theta}
\sum^{N+1}_{z=0} \binom{N+1}{z} \theta^z (1-\theta)^{ N+1- z}
\sum_{j=0}^{2l} \binom{2l}{j} \{z - (N+1)\theta \}^j \{- (1-\theta)\}^{2l-j} \\
&\le \frac{\displaystyle \max \left( 1, \; \frac{1}{a} \right)}{(N \theta +a)^{2l-1} (N+1)\theta}
\sum_{j=0}^{2l} \binom{2l}{j} |\mu_j(N+1,\theta)|,
\end{align*}
where we define $\mu_0(N+1,\theta):=1$.
By Lemma 2, there exist positive constants $\tilde{C}_i$ $(i=0,\ldots,2l)$ such that
\begin{align*}
 \mathrm{E} \biggl( & \frac{w^{2l}} {1 + w} \biggr)
\le \frac{\displaystyle \max \left( 1, \; \frac{1}{a} \right)}{(N \theta +a)^{2l-1} (N+1)\theta}
\left\{
\sum_{j=0}^{l} \tilde{C}_{2j} \{(N+1) \theta\}^j
+
\sum_{j=0}^{l-1} \tilde{C}_{2j+1} \{(N+1) \theta\}^j
\right\} \\
\le& \max \left( 1, \; \frac{1}{a} \right)
\left\{
\sum_{j=0}^{l} \left(\frac{N+1}{N}\right)^{j-1} \frac{\tilde{C}_{2j}}{(N \theta)^{2l-j}}
+
\sum_{j=0}^{l-1} \left(\frac{N+1}{N}\right)^{j-1} \frac{\tilde{C}_{2j+1}}{(N \theta)^{2l-j}}
\right\} \\
\le& \; 2^{l-1} \max \left( 1, \; \frac{1}{a} \right) \frac{1}{(N \theta)^{l}}
\left\{
\tilde{C}_{2l} +
\sum_{j=0}^{l-1} \frac{\tilde{C}_{2j}}{(N \epsN)^{l-j}}
+
\sum_{j=0}^{l-1} \frac{\tilde{C}_{2j+1}}{(N \epsN)^{l-j}}.
\right\}
\end{align*}
Since $N \epsN$ goes to infinity, there exists a constant $C$ such that
\begin{align*}
\mathrm{E} & \left(  \frac{w^{2l}} { 1 + w} \right)
\le \frac{C}{(N\theta)^l}.
\end{align*}
Therefore,
\[
\mathrm{E} \left( - \frac{ w^{2l+1}}{ 1+w} \right)
= - \E(w^{2l}) + \E \left(\frac{w^{2l}}{1 + w}\right)
\le \frac{C} {(N\theta)^l}.
\] 
\qed

\vspace{0.5cm}

\noindent
Proof of Lemma 4.
The desired inequality is equivalent to
\begin{align}
(\alpha_1+\cdots+\alpha_k) B_{\Delta_\eps} (\alpha_1 + 1, \alpha_2, \dotsb, \alpha_k) - \alpha_1 B_{\Delta_\eps} (\alpha_1, \dotsb, \alpha_k)
\leq& \;
\frac{\Gamma(\alpha_2) \cdots \Gamma(\alpha_k)}{\Gamma(\alpha_2+\cdots+\alpha_k)}
\varepsilon ^{\alpha_1} (1-\varepsilon)^{\alpha_2+\dotsb+\alpha_k}.
\label{lm4equiv}
\end{align}

Let
\[
 w_i := \frac{\theta_i}{1 - \theta_1} ~~~~~ (i=2,\ldots,k).
\]
Then,
\begin{align*}
\theta_1^{\alpha_1-1} & \cdots \theta_{k-1}^{\alpha_{k-1}-1} (1-\theta_1-\cdots-\theta_{k-1})^{\alpha_k-1}
\dd \theta_1 \cdots \dd \theta_{k-1} \\
=& \theta_1^{\alpha_1-1} (1-\theta_1)^{\alpha_2 + \cdots + \alpha_k - 1} 
w_2^{\alpha_2-1} \cdots w_{k-1}^{\alpha_{k-1}-1} (1-w_2-\cdots-w_{k-1})^{\alpha_k-1}
\dd \theta_1 \dd w_2 \cdots \dd w_{k-1}.
\end{align*}

We define
\[
 {\Delta'_{\eps/(1-\eps)}} := \{ w = (w_2,\ldots,w_{k-1}) \mid w_i \geq \eps/(1-\eps) ~\; (i=2,\ldots,k),~ w_k := 1-\sum_{i=2}^{k-1} w_i \}.
\]
If $\theta \in \Delta_\eps$, then $(\theta_2/(1-\theta_1),\ldots,\theta_{k-1}/(1-\theta_1)) \in {\Delta'_{\eps/(1-\eps)}}$.

If $w = (w_2,\ldots,w_{k-1}) \in \Delta'_{\varepsilon/(1-\varepsilon)}$ is fixed,
then
$\{\theta = (\theta_1,\ldots,\theta_{k-1}) \mid \theta \in \Delta_\eps, ~ \theta_i = (1-\theta_1) w_i ~~(i=2,\ldots,k-1)\}$,
which is a subset of $\Delta_\eps$,
is represented as
$\{\theta \mid (\theta_1,(1-\theta_1)w_2,\ldots,(1-\theta_1)w_{k-1}) \mid L(w_2,\ldots,w_{k-1}) \leq \theta_1 \leq U(w_2,\ldots,w_{k-1}) \}$
by using appropriate functions $L(w_2,\ldots,w_{k-1})$ and $U(w_2,\ldots,w_{k-1})$
because $\Delta_\eps$ is a bounded closed convex set.

If $(\theta_1, (1-\theta_1) w_2, \ldots, (1-\theta_1) w_k) \in \Delta_\eps$,
then $(\eps, (1-\eps) w_2, \ldots, (1-\eps) w_{k-1}) \in \Delta_\eps$
because $(1-\eps) w_i \geq (1-\theta_1) w_i \geq \eps$ for $i=2,\ldots,k$
and $\eps + \sum_{i=2}^k (1-\eps) w_i = 1$.
Thus, $L(w_2,\ldots,w_{k-1}) \geq \eps$.
Obviously, $L(w_2,\ldots,w_{k-1}) \le \eps$ because $\theta \notin \Delta_\eps$ if $\theta_1 < \eps$.
Hence, $L(w_2,\ldots,w_{k-1}) = \eps$.

Since $\theta_1 = 1 - \theta_2 - \cdots - \theta_k \leq 1 - (k-1) \eps$,
$U(w_2,\ldots,w_{k-1}) \leq 1- (k-1) \eps$.

Therefore, we have
\begin{align*}
B_{\Delta_\eps} & (\alpha_1 + 1, \alpha_2, \dotsb, \alpha_k)
= \int_{\Delta_\eps}
\theta_1^{\alpha_1} \theta_2^{\alpha_2-1} \cdots \theta_{k-1}^{\alpha_{k-1}-1} (1-\theta_1-\cdots-\theta_{k-1})^{\alpha_k-1}
\dd \theta_1 \cdots \dd \theta_{k-1} \\
=& \int_{\Delta'_{\eps/(1-\eps)}}
\left\{ \int_\varepsilon^{U(w_2,\ldots,w_{k-1})}
\hspace{-20pt} \theta_1^{\alpha_1} (1-\theta_1)^{\sum_{i=2}^k \alpha_i - 1} \dd \theta_1 \right\} \\
& \times w_2^{\alpha_2-1} \cdots w_{k-1}^{\alpha_{k-1}-1} (1-w_2-\cdots-w_{k-1})^{\alpha_k-1}
\dd w_2 \cdots \dd w_{k-1} \\
=& \int_{\Delta'_{\eps/(1-\eps)}}
\left\{ \left[ \theta_1^{\alpha_1} \frac{-1}{\sum_{i=2}^k \alpha_i}
(1-\theta_1)^{\sum_{i=2}^k \alpha_i} \right]^{U(w_2,\ldots,w_{k-1})}_\varepsilon
\hspace{-15pt} + \int_\varepsilon^{U(w_2,\ldots,w_{k-1})} \hspace{-10pt}
\frac{\alpha_1}{\sum_{i=2}^k \alpha_i} \theta_1^{\alpha_1-1} (1-\theta_1)^{\sum_{i=2}^k \alpha_i} \dd \theta_1 \right\} \\
& \times w_2^{\alpha_2-1} \cdots w_{k-1}^{\alpha_{k-1}-1} (1-w_2-\cdots-w_{k-1})^{\alpha_k-1}
\dd w_2 \cdots \dd w_{k-1} \\
=& \int_{\Delta'_{\eps/(1-\eps)}}
\frac{1}{\sum_{i=2}^k \alpha_i} \biggl[
\varepsilon^{\alpha_1} (1-\varepsilon)^{\sum_{i=2}^k \alpha_i}
- \{U(w_2,\ldots,w_{k-1})\}^{\alpha_1} \{1-U(w_2,\ldots,w_{k-1})\}^{\sum_{i=2}^k \alpha_i}
\biggr] \\
& \times w_2^{\alpha_2-1} \cdots w_{k-1}^{\alpha_{k-1}-1} (1-w_2-\cdots-w_{k-1})^{\alpha_k-1}
\dd w_2 \cdots \dd w_{k-1} \\
& + \int_{\Delta'_{\eps/(1-\eps)}}
\left\{ \int_\varepsilon^{U(w_2,\ldots,w_{k-1})} 
\vspace{-20pt} \frac{\alpha_1}{\sum_{i=2}^k \alpha_i} \theta_1^{\alpha_1-1} (1-\theta_1)^{\sum_{i=2}^k \alpha_i} \dd \theta_1 \right\} \\
& \times w_2^{\alpha_2-1} \cdots w_{k-1}^{\alpha_{k-1}-1} (1-w_2-\cdots-w_{k-1})^{\alpha_k-1}
\dd w_2 \cdots \dd w_{k-1} \\
\leq& \frac{1}{\sum_{i=2}^k \alpha_i}
\varepsilon^{\alpha_1} (1-\varepsilon)^{\sum_{i=2}^k \alpha_i} 
\int_{\Delta'_{\eps/(1-\eps)}} w_2^{\alpha_2-1} \cdots w_{k-1}^{\alpha_{k-1}-1} (1-w_2-\cdots-w_{k-1})^{\alpha_k-1}
\dd w_2 \cdots \dd w_{k-1} \\
& + \frac{\alpha_1}{\sum_{i=2}^k \alpha_i} B_{\Delta'_{\eps/(1-\eps)}} (a_1, a_2, \dotsb, a_k) \\
\leq& \; \frac{B(\alpha_2,\ldots,\alpha_k)}{\sum_{i=2}^k \alpha_i}
\varepsilon^{\alpha_1} (1-\varepsilon)^{\alpha_2+\cdots+\alpha_k}
+ \frac{\alpha_1}{\sum_{i=2}^k \alpha_i} B_{\Delta_\eps} (a_1, a_2, \dotsb, a_k).
\end{align*}
Thus, \eqref{lm4equiv} is obtained.
\qed

\vspace{0.5cm}

\noindent
Proof of Lemma 5.
We obtain the desired inequality from
\begin{align*}
\frac{\partial}{\partial t} \frac{B_{[s,t]}(\alpha+1,\beta)}{B_{[s,t]}(\alpha,\beta)}
=& \frac{\displaystyle \left\{\frac{\partial}{\partial t} B_{[s,t]}(\alpha+1,\beta) \right\} B_{[s,t]}(\alpha,\beta)
- B_{[s,t]}(\alpha+1,\beta) \left\{\frac{\partial}{\partial t} B_{[s,t]}(\alpha,\beta) \right\}}
{\displaystyle \{B_{[s,t]}(\alpha,\beta)\}^2} \\
=& \frac{1}{\displaystyle \{B_{[s,t]}(\alpha,\beta)\}^2}
\biggl\{ t^\alpha (1 - t)^{\beta-1} B_{[s,t]}(\alpha,\beta)
- t^{\alpha-1} (1 - t)^{\beta-1} B_{[s,t]}(\alpha+1,\beta) \biggr\} \\
=& \frac{1}{\displaystyle \{B_{[s,t]}(\alpha,\beta)\}^2}
t^{\alpha-1} (1- t)^{\beta-1}
\int_s^t (t - \theta) \theta^{\alpha-1} (1-\theta)^{\beta-1} \dd \theta \geq 0
\end{align*}
and
\begin{align*}
\frac{\partial}{\partial s} \frac{B_{[s,t]}(\alpha+1,\beta)}{B_{[s,t]}(\alpha,\beta)}
=& \frac{\displaystyle \left\{\frac{\partial}{\partial s} B_{[s,t]}(\alpha+1,\beta) \right\} B_{[s,t]}(\alpha,\beta)
- B_{[s,t]}(\alpha+1,\beta) \left\{\frac{\partial}{\partial s} B_{[s,t]}(\alpha,\beta) \right\}}
{\displaystyle \{B_{[s,t]}(\alpha,\beta)\}^2} \\
=& \frac{1}{\displaystyle \{B_{[s,t]}(\alpha,\beta)\}^2}
\biggl\{ - s^a (1-s)^{\beta-1} B_{[s,t]}(\alpha,\beta)
- \{-s^{\alpha-1} (1-s)^{\beta-1} \bigr\} B_{[s,t]}(\alpha+1,\beta) \biggr\} \\
=& \frac{1}{\displaystyle \{B_{[s,t]}(\alpha,\beta)\}^2}
s^{\alpha-1} (1-s)^{\beta-1}
\int_s^t (\theta - s) \theta^{\alpha-1} (1-\theta)^{\beta-1} \dd \theta
\geq 0.
\end{align*}
\qed

\vspace{0.5cm}

\noindent
Proof of Lemma 6.
Define $w_i$ $(i=2,\ldots,k)$, ${\Delta'_{\eps/(1-\eps)}}$, and $U(w_2,\ldots,w_{k-1})$ as in the proof of Lemma 4.
Let
\begin{align*}
\tilde{p}(\theta_1,&w_2,\ldots,w_{k-1}) \dd \theta_1 \dd w_2 \cdots \dd w_{k-1}
:= \theta_1^{\alpha_1-1} \cdots \theta_{k-1}^{\alpha_{k-1}-1} (1-\theta_1-\cdots-\theta_{k-1})^{\alpha_k-1}
\dd \theta_1 \cdots \dd \theta_{k-1} \\
=& \theta_1^{\alpha_1-1} (1-\theta_1)^{\alpha_2 + \cdots + \alpha_k - 1} 
w_2^{\alpha_2-1} \cdots w_{k-1}^{\alpha_{k-1}-1} (1-w_2-\cdots-w_{k-1})^{\alpha_k-1}
\dd \theta_1 \dd w_2 \cdots \dd w_{k-1}.
\end{align*}
and
\begin{align*}
p(\theta_1,w_2,\ldots,w_{k-1}) :=&
\frac{\tilde{p}(\theta_1,w_2,\ldots,w_{k-1})}{B_{\Delta_\varepsilon}(\alpha_1,\alpha_2,\ldots,\alpha_k)}.
\end{align*}
Since
\begin{align*}
B_{\Delta_\varepsilon}(\alpha_1,\alpha_2,\ldots,\alpha_k)
:= \int_{\Delta'_{\eps/(1-\eps)}}
\left\{\int_\varepsilon^{U(w_2,\ldots,w_{k-1})} \tilde{p}(\theta_1,w_2,\ldots,w_{k-1}) \dd \theta_1 \right\} \dd w_2 \cdots \dd w_{k-1},
\end{align*}
$p(\theta_1,w_2,\ldots,w_{k-1})$ is a probability density.
The marginal density of $(w_2,\ldots,w_{k-1})$ is
\begin{align*}
p(w_2,&\ldots,w_k) = \int_\varepsilon^{U(w_2,\ldots,w_k)} p(\theta_1,w_2,\ldots,w_k) \dd \theta_1.
\end{align*}
The conditional density of $\theta_1$ given $(w_2,\ldots,w_{k-1})$ is
\begin{align*}
p(\theta_1 & \mid w_2,\ldots,w_{k-1}) = \frac{p(\theta_1,w_2,\ldots,w_{k-1})}{p(w_2,\ldots,w_{k-1})} \\
=& 
\begin{cases}
\displaystyle \frac{\theta_1^{\alpha_1-1}(1-\theta_1)^{\alpha_2+\cdots+\alpha_{k}-1}}
{\displaystyle
\int_\varepsilon^{U(w_2,\ldots,w_{k-1})} \theta_1^{\alpha_1-1}(1-\theta_1)^{\alpha_2+\cdots+\alpha_{k}-1} \dd \theta_1},
 & \varepsilon \leq \theta_1 \leq U(w_2,\ldots,w_{k-1}), \\
0, & \text{otherwise}.
\end{cases}
\end{align*}
Then, from Lemma 5 and $U(w_2,\ldots,w_{k-1}) \leq 1-(k-1) \eps \leq 1$,
\begin{align*}
\lefteqn{\frac{B_{\Delta_\varepsilon}(\alpha_1+1,\alpha_2,\ldots,\alpha_k)}
{B_{\Delta_\varepsilon}(\alpha_1,\alpha_2,\ldots,\alpha_k)} =
\int_{\Delta'_{\varepsilon/(1-\varepsilon)}}
\int_{\varepsilon}^{U(w_2,\ldots,w_{k-1})} \theta_1
p(\theta_1,w_2,\ldots,w_{k-1}) \dd \theta_1 \dd w_2 \cdots \dd w_{k-1}} \\
&= \int_{\Delta'_{\eps/(1-\eps)}} \left\{
\int_{\varepsilon}^{U(w_2,\ldots,w_{k-1})} \theta_1
p(\theta_1 \mid w_2,\ldots,w_{k-1}) \dd \theta_1 \right\}
p(w_2,\ldots,w_{k-1}) \dd w_2 \cdots \dd w_{k-1} \\
&= \int_{\Delta'_{\eps/(1-\eps)}} \left\{
\frac{\displaystyle
\int_{\varepsilon}^{U(w_2,\ldots,w_{k-1})} \theta_1^{\alpha_1} (1-\theta_1)^{\alpha_2+\cdots+\alpha_{k}-1} \dd \theta_1}
{\displaystyle
\int_\varepsilon^{U(w_2,\ldots,w_{k-1})} \theta_1^{\alpha_1-1}(1-\theta_1)^{\alpha_2+\cdots+\alpha_{k}-1} \dd \theta_1}
\right\}
p(w_2,\ldots,w_{k-1}) \dd w_2 \cdots \dd w_{k-1} \\
&\leq
\frac{\displaystyle \int_\varepsilon^{1-(k-1)\varepsilon}
{\theta_1}^{\alpha_1} (1-\theta_1)^{\alpha_2+\cdots+\alpha_k-1} \dd \theta_1}
{\displaystyle \int_\varepsilon^{1-(k-1)\varepsilon}
{\theta_1}^{\alpha_1-1} (1-\theta_1)^{\alpha_2+\cdots+\alpha_k-1} \dd \theta_1}
=
\frac{B_{[\varepsilon,1-(k-1)\varepsilon]}(\alpha_1+1,\alpha_2+\cdots+\alpha_k)}
{B_{[\varepsilon,1-(k-1)\varepsilon]}(\alpha_1,\alpha_2+\cdots+\alpha_k)} \\
&\le \frac{B_{[\varepsilon,1]}(\alpha_1+1,\alpha_2+\cdots+\alpha_k)}{B_{[\varepsilon,1]}(\alpha_1,\alpha_2+\cdots+\alpha_k)}.
\end{align*}
\begin{flushright}
\vspace{-15pt} 

\qed
\end{flushright}

\vspace{0.5cm}

\noindent
Proof of Lemma 7.
The right hand side of the equation is represented by
\begin{align*}
& \hspace{-30pt} \sum_{x_2,\ldots,x_{k}: \sum_{i=2}^{k} x_i = N-x_1}
\frac{B(x_1+\alpha_1,x_2+\alpha_2,\ldots,x_k+\alpha_k)}{B(\alpha_1,\alpha_2,\ldots,\alpha_k)} \binom{N}{x_1, \ldots, x_k} \\
= & \sum_{x_2,\ldots,x_{k}: \sum_{i=2}^{k} x_i = N-x_1}
\int_\Delta \frac{\theta_1^{\alpha_1-1} \cdots \theta_{k-1}^{\alpha_{k-1}-1} (1-\sum_{i=1}^{k-1} \theta_i)^{\alpha_k-1}}
{B(\alpha_1,\ldots,\alpha_k)} \\
& \times \binom{N}{x_1,\ldots,x_k} \theta_1^{x_1} \cdots \theta_{k-1}^{x_{k-1}} (1-\sum_{i=1}^{k-1} \theta_i)^{N-\sum_{i=1}^{k-1} x_i}
\dd \theta_1 \cdots \dd \theta_{k-1} \\
= &
\int_\Delta \frac{\theta_1^{\alpha_1-1} \cdots \theta_{k-1}^{\alpha_{k-1}-1} (1-\sum_{i=1}^{k-1} \theta_i)^{\alpha_k-1}}{B(\alpha_1,\ldots,\alpha_k)}
\binom{N}{x_1} \theta_1^{x_1} (1-\theta_1)^{N-x_1}
\dd \theta_1 \cdots \dd \theta_{k-1}.
\end{align*}
From the relation
\begin{align*}
\int_\Delta & \theta_1^{\alpha_1-1} \cdots \theta_{k-1}^{\alpha_{k-1}-1} (1-\sum_{i=1}^{k-1} \theta_i)^{\alpha_k-1}
\theta_1^{x_1} (1-\theta_1)^{N-x_1}
\dd \theta_1 \cdots \dd \theta_{k-1} \\
=&
\int_0^1 \int_0^{1-\theta_1} \cdots \int_0^{1-\sum_{i=1}^{k-2} \theta_i}
\theta_1^{\alpha_1-1} \cdots \theta_{k-1}^{\alpha_{k-1}-1} (1-\sum_{i=1}^{k-1} \theta_i)^{\alpha_k-1}
\theta_1^{x_1} (1-\theta_1)^{N-x_1} \dd \theta_{k-1} \dd \theta_{k-2} \cdots \dd \theta_1 \\
=& \int_0^1 \int_0^{1-\theta_1} \cdots \int_0^{1-\sum_{i=1}^{k-3} \theta_i} \int_0^{1}
 \theta_1^{\alpha_1 + x_1 -1} \theta_2^{\alpha_2-1} \cdots \theta_{k-2}^{\alpha_{k-2}-1} (1-\sum_{i=1}^{k-2} \theta_i)^{\alpha_{k-1}+\alpha_k-1}
\btheta_{k-1}^{\alpha_{k-1}-1} (1- \btheta_{k-1})^{\alpha_k-1} \\
& \times (1-\theta_1)^{N-x_1} \dd \btheta_{k-1} \dd \theta_{k-2} \cdots \dd \theta_1 \\
=& B(\alpha_{k-1}, \alpha_k) \int_0^1 \int_0^{1-\theta_1} \cdots \int_0^{1-\sum_{i=1}^{k-3} \theta_i}
 \theta_1^{\alpha_1 + x_1 -1} \theta_2^{\alpha_2-1} \cdots \theta_{k-2}^{\alpha_{k-2}-1} (1-\sum_{i=1}^{k-2} \theta_i)^{\alpha_{k-1}+\alpha_k-1} \\
& \times (1-\theta_1)^{N-x_1} \dd \theta_{k-2} \cdots \dd \theta_1 \\
=& B(\alpha_{k-1}, \alpha_k) B(\alpha_{k-2}, \alpha_{k-1}+\alpha_k) \cdots B(\alpha_2, \sum_{i=3}^k \alpha_i)
 \int_0^1 \theta_1^{\alpha_1 + x_1 -1} (1 - \theta_1)^{\sum_{i=2}^{k} \alpha_i + N - x_1 - 1} \dd \theta_1 \\
=& B(\alpha_2, \alpha_3, \ldots, \alpha_k) B(\alpha_1+x_1,\sum_{i=2}^{k} \alpha_i + N - x_1) \\
=& \frac{B(\alpha_1, \alpha_2, \ldots, \alpha_k)}{B(\alpha_1, \sum_{i=2}^k \alpha_i)} B(\alpha_1+x_1,\sum_{i=2}^{k} \alpha_i + N - x_1),
\end{align*}
where $\btheta_{k-1} = \theta_{k-1}/ \sum_{j=1}^{k-2} \theta_j$,
we obtain the desired result.
\qed

\vspace{0.5cm}

\noindent
Proof of Lemma 8.
Since
\begin{align*}
\int^1_\varepsilon & \theta^\alpha (1-\theta)^{\beta-1} \dd \theta
= \left[-\frac{1}{\beta} \theta^\alpha (1-\theta)^\beta \right]^1_\varepsilon
+ \frac{\alpha}{\beta} \int^1_\varepsilon \theta^{\alpha-1} (1-\theta)^\beta \dd \theta \\
=& \frac{1}{\beta} \varepsilon^\alpha (1-\varepsilon)^\beta
- \frac{\alpha}{\beta} \int^1_\varepsilon \theta^\alpha (1-\theta)^{\beta-1} \dd \theta
+ \frac{\alpha}{\beta} \int^1_\varepsilon \theta^{\alpha-1} (1-\theta)^{\beta-1} \dd \theta,
\end{align*}
we have
\begin{align*}
(\alpha+\beta) & \frac{\int^1_\varepsilon \theta^\alpha (1-\theta)^{\beta-1} \dd \theta}{\int^1_\varepsilon \theta^{\alpha-1} (1-\theta)^{\beta-1} \dd \theta}
= \alpha + \frac{\varepsilon^\alpha (1-\varepsilon)^\beta}{\int^1_\varepsilon \theta^{\alpha-1} (1-\theta)^{\beta-1} \dd \theta}
\le \alpha + \frac{\varepsilon^\alpha (1-\varepsilon)^\beta}{\varepsilon^{\alpha-1} \int^1_\varepsilon (1-\theta)^{\beta-1} \dd \theta} \\
=& \alpha + \frac{\varepsilon (1-\varepsilon)^\beta}{\frac{-1}{\beta} [ (1-\theta)^\beta]^1_\varepsilon}
= \alpha + \beta \varepsilon.
\end{align*}
\qed
\vspace{0.5cm}

\begin{center}
Acknowledgments
\end{center}
This research was partially supported
by Grant-in-Aid for Scientific Research (23300104)
and by the Aihara Project, the FIRST program from JSPS, initiated by CSTP.

\bibliographystyle{plainnat}

\end{appendix}
\end{document}